# CAPACITY ON FINSLER SPACES[*]

## B. BIDABAD[1][**] AND S. HEDAYATIAN[2]

[1]Faculty of Mathematics and Computer Sciences, Amirkabir University of Technology,
424 Hafez Ave, 15914 Tehran, I. R. of Iran
Email: bidabad@aut.ac.ir

[1]Faculty of Mathematics and Computer Sciences, Chamran University of Ahvaz, Golestan Bld,
P.O. Box: 61355-83151 Ahvaz, I. R. of Iran
Email: Hedayatian@scu.ac.ir

**Abstract** – Here, the concept of electric capacity on Finsler spaces is introduced and the fundamental conformal invariant property is proved, i.e. the capacity of a compact set on a connected non-compact Finsler manifold is conformal invariant. This work enables mathematicians and theoretical physicists to become more familiar with the global Finsler geometry and one of its new applications.

**Keywords** – Capacity, conformal invariant, Finsler space

## 1. INTRODUCTION

Finsler space is the most natural and advanced generalization of Euclidean space, which has many applications in theoretical physics. The physical notion of capacity is the electrical capacity of a $2$-dimensional conducting surface, which is defined as the ratio of a given positive charge on the conductor to the value of the potential on its surface.

The capacity of a set as a mathematical concept was introduced first by N. Wiener in 1924 and was subsequently developed by O. Forstman [1], C. J. de La Vallee Poussin, and several other physicists and mathematicians in connection with the potential theory.

The concept of conformal capacity was introduced by Loewner [2] and has been extensively developed for $\mathbb{R}^n$ [3-6]. In particular, it was used by G.D. Mostow to prove his famous theorem on the rigidity of hyperbolic spaces [5]. The concept of capacity on Riemannian geometry was introduced by J. Ferrand [7] and developed in the joint work's of M. Vuorinan and G.J. Martin [8] and [9].

Here, we introduce the concept of capacity for Finsler spaces and prove that, it depends only on the conformal structure of $(M, g)$, more precisely:

**Theorem:** Let $(M, g)$ be a connected non-compact Finsler manifold, then the capacity of a compact set on $M$ is a conformal invariant.

# 1. PRELIMINARIES

## 1.1. Finsler metric

Let $M$ be an n-dimensional $C^\infty$ manifold. For a point $x \in M$, denote by $T_xM$ the tangent space of $M$ at $x$. The tangent bundle $TM$ on $M$ is the union of tangent spaces $T_xM$. We will denote the elements of TM by $(x, y)$ where $y \in T_xM$. Let $TM_0 = TM \setminus \{0\}$. The natural projection $\pi : TM \to M$ is given by $\pi(x, y) := x$. Throughout this paper we use the *Einstein summation convention* for the expressions with repeated indices. That is, wherever an index appears twice, once as a subscript, and once as a superscript, then that term is summed over all values of that index.

A *Finsler structure* on a manifold $M$ is a function $F : TM_0 \to [0, \infty)$ with the following properties: (i) $F$ is $C^\infty$ on $TM_0$. (ii) $F$ is positively 1-homogeneous on the fibers of tangent bundle $TM$, i.e. $\forall \lambda > 0 \quad F(x, \lambda y) = \lambda F(x, y)$. (iii) The Hessian of $F^2$ with elements $g_{ij}(x, y) := \frac{1}{2}[F^2(x, y)]_{y^i y^j}$ is positive definite on $TM_0$. We recall that, $g_{ij}$ is a homogeneous tensor of degree zero in $y$ and $g_{ij}(x, y)y^i y^j = g(y, y)$, where $g( , )$ is the local scalar product on any point of $TM_0$. Then the pair $(M, g)$ is called a *Finsler manifold*. The Finsler structure $F$ is Riemannian if $g_{ij}(x, y)$ are independent of $y \neq 0$.

## 1.2. Notations on conformal geometry of Finsler manifolds

Let's consider two $n$-dimensional Finsler manifolds $(M, g)$ and $(M', g')$ with Finsler structures $F$ and $F'$ and with line elements $(x, y)$ and $(x', y')$ respectively. Throughout this paper we shall assume that coordinate systems on $(M, g)$ and $(M', g')$ have been chosen so that $x'^i = x^i$ and $y'^i = y^i$ holds for all $i$, unless a contrary assumption is explicitly made. Using this assumption these manifolds can be denoted simply by $M$ and $M'$, respectively. Let $u$ and $v$ be two tangent vectors at a point $x$ of a Finsler manifold $(M, g)$. The *angle $\theta$ of $v$ with respect to $u$* is defined by

$$\cos \theta = \frac{g_{ij}(x, u) u^i v^j}{\sqrt{g_{ij}(x, u) u^i u^j} \sqrt{g_{ij}(x, u) v^i v^j}}.$$

Clearly this notion of angle is not symmetric. A diffeomorphism $f : M \to M'$ between two Finsler manifolds is called *conformal* if for each $p \in M$, $(f_*)_p$ preserves the angles of any tangent vector, with respect to any $y$ in $M$. In this case the two Finsler manifolds are called *conformal equivalent* or simply *conformal*. If $M = M'$ then $f$ is called a *conformal transformation* or *conformal automorphism*. It can be easily checked that a diffeomorphism is conformal if and only if $f^* g' = e^{2\sigma} g$ for some function $\sigma : M \to IR$ (this result is due to Knebelman [10]. In fact, the sufficient condition implies that the function $\sigma(x, y)$ be independent of direction $y$, or equivalently $\partial \sigma / \partial y^i = 0$). The diffeomorphism $f$ is called an *isometry* if $f^* g' = g$. Two Finsler structures $F$ and $F'$ are called *conformal* if $F'(x, y) = e^\sigma F(x, y)$ or equivalently, $g' = e^{2\sigma(x)} g$. Locally we have $g'_{ij}(x, y) = e^{2\sigma(x)} g_{ij}(x, y)$, and $g'^{ij}(x, y) = e^{-2\sigma(x)} g^{ij}(x, y)$.

## 1.3. Some vector bundles and their properties

Let $\pi : TM \longrightarrow M$ be the natural projection from $TM$ to $M$. The *pull-back tangent space* $\pi^* TM$ is defined by $\pi^* TM := \{(x, y, v) \mid y \in T_xM_0, v \in T_xM\}$. The *pull-back cotangent space* $\pi^* T^* M$ is the dual of $\pi^* TM$. Both $\pi^* TM$ and $\pi^* T^* M$ are n-dimensional vector spaces over $TM_0$ [11, 12]. We denote by $S_xM$ the set consisting of all rays $[y] := \{\lambda y \mid \lambda > 0\}$, where $y \in T_xM_0$. Let $SM = \bigcup_{x \in M} S_xM$, then $SM$ has a natural $(2n-1)$ dimensional manifold structure and the total space of a fiber bundle, called *Sphere bundle* over $M$. We denote the elements of $SM$ by $(x, [y])$ where $y \in T_x M_0$. Let $p: SM \longrightarrow M$ denote the natural projection from $SM$ to $M$. The *pull-back tangent space* $p^*TM$ is defined by $p^*TM := \{(x,[y],v) | y \in T_x M_0, v \in T_x M\}$. The *pull-back cotangent space* $p^*T^*M$ is the dual of $p^*TM$. Both $p^*TM$ and $p^*T^*M$ are total spaces of vector bundles over $SM$. We use the following Lemma for replacing the $C^\infty$ functions on $TM_0$ by those on $SM$.

**Lemma 1.1.** [13] Let $\eta$ be the function $\eta: TM_0 \longrightarrow SM$, where $\eta(x,y) = (x,[y])$ and $f \in C^\infty(TM_0)$. Then there exists a function $g \in C^\infty(SM)$ satisfying $\eta^*g = f$ if and only if $f(x,y) = f(x,\lambda y)$, where $y \in T_x M_0, \lambda > 0$ and $\eta^*$ is the pull-back of $\eta$.

Let $f \in C^\infty(M)$, the vertical lift of $f$ denoted by $f^V \in C^\infty(TM_0)$, be defined by $f^V: TM \longrightarrow IR$, where $f^V(x,y) := f \circ \pi(x,y) = f(x)$. $f^V$ is independent of $y$ and from Lemma 1.1 there is a function $g$ on $C^\infty(SM)$ related to $f^V$ by means of $\eta^*g = f^V$. In the sequel $g$ is denoted by $f^V$ for simplicity. It is well known that, if the differentiable manifold $M$ is compact then the Sphere bundle $SM$ is compact, and also it is orientable whether $M$ is orientable or not [14, 15]).

### 1.4. Nonlinear connections

#### 1.4.1. Nonlinear connection on the tangent bundle TM

Consider $\pi_*: TTM \longrightarrow TM$ and put $\ker \pi_*^v = \{z \in TTM | \pi_*^v(z) = 0\}$, $\forall v \in TM$, then the vertical vector bundle on $M$ is defined by $VTM = \bigcup_{v \in TM} \ker \pi_*^v$. A *non-linear connection* or a *horizontal distribution* on $TM$ is a complementary distribution $HTM$ for $VTM$ on $TTM$. These functions are called coefficients of the non-linear connection and will be noted in the sequel by $N_i^j$. It is clear that $HTM$ is a vector sub-bundle of $TTM$ called horizontal vector bundle. Therefore we have the decomposition $TTM = VTM \oplus HTM$.

Using the induced coordinates $(x^i, y^i)$ on $TM$, where $x^i$ and $y^i$ are called, respectively, *position* and *direction* of a point on $TM$, we have the local field of frames $\{\frac{\partial}{\partial x_i}, \frac{\partial}{\partial y_i}\}$ on $TTM$. Let $\{dx^i, dy^i\}$ be the dual of $\{\frac{\partial}{\partial x^i}, \frac{\partial}{\partial y^i}\}$. It is well known that we can choose a local field of frames $\{\frac{\delta}{\delta x^i}, \frac{\partial}{\partial y_i}\}$ adapted to the above decomposition, i.e. $\frac{\delta}{\delta x^i} \in \chi(HTM)$ and $\frac{\partial}{\partial y_i} \in \chi(VTM)$. They are sections of horizontal and vertical bundles, $HTM$ and $VTM$, defined by $\frac{\delta}{\delta x^i} = \frac{\partial}{\partial x_i} - N_i^j \frac{\partial}{\partial y_j}$, where $N_i^j(x,y)$ are the coefficients of non linear $\gamma^i_{jk} := \frac{1}{2} g^{is}(\frac{\partial g_{sj}}{\partial x^k} - \frac{\partial g_{jk}}{\partial x^s} + \frac{\partial g_{ks}}{\partial x^j})$ and $C_{ijk} = \frac{1}{2} \frac{\partial g_{ij}}{\partial y^k}$.

#### 1.4.2. Nonlinear connections on the sphere bundle SM

Using the coefficients of non linear connection on $TM$, one can define a non linear connection on $SM$ by using the objects which are invariant under positive re-scaling $y \mapsto \lambda y$. Our preference for remaining on $SM$ forces us to work with $\frac{N^i_j}{F} := \gamma^i_{jk} l^k - C^i_{jk} \gamma^k_{rs} l^r l^s$, where $l^i = \frac{y^i}{F}$. We also prefer to work with the local field of frames $\{\frac{\delta}{\delta x^i}, F\frac{\partial}{\partial y^j}\}$ and $\{dx^i, \frac{\delta y^j}{F}\}$, which are invariant under the positive re-scaling of $y$, and therefore, live over $SM$. They can also be used as a local field of frames over tangent bundle $p^*TM$ and cotangent bundle $p^*T^*M$ respectively.

### 1.5. A Riemannian metric on SM

It turns out that the manifold $TM_0$ has a natural Riemannian metric, known in the literature as *Sasaki metric* [12, 16]); $\tilde{g} = g_{ij}(x,y) dx^i \otimes dx^j + g_{ij}(x,y) \frac{\delta y^i}{F} \otimes \frac{\delta y^j}{F}$, where $g_{ij}(x,y)$ is the Hessian of Finsler structure $F^2$. They are functions on $TM_0$ and invariant under positive re-scaling of $y$, therefore they can be considered as functions on $SM$. With respect to this metric, the *horizontal subspace* spanned by $\frac{\delta}{\delta x^j}$ is orthogonal to the *vertical subspace* spanned by $F\frac{\partial}{\partial y^i}$. The metric $\tilde{g}$ is invariant under the positive re-scaling of $y$ and can be considered as a Riemannian metric on $SM$.

### 1.6. Hilbert form

Consider the pull-back vector bundle $p^*TM$ over $SM$. The pull-back tangent bundle $p^*TM$ has a canonical section $l$ defined by $l_{(x,[y])} = (x,[y],\frac{y}{F(x,y)})$. We use the local coordinate system $(x^i, y^i)$ for $SM$, where $y^i$ are homogeneous coordinates up to a positive factor. Let $\{\partial_i\}$ be a natural local field of frames for $p^*TM$, where $\partial_i := (x,[y],\frac{\partial}{\partial x^i})$. The natural dual co-frame for $p^*T^*M$ is noted by $\{dx^i\}$. The Finsler structure $F(x,y)$ induces a canonical 1-form on $SM$ defined by $\omega := l_i dx^i$, where $l_i = g_{ij}l^j$ and $\omega$ is called the *Hilbert form* of $F$. Using $g_{ij} = FF_{y^i y^j} + F_{y^i} F_{y^j}$ and $\frac{\delta F}{\delta x^i} = 0$, with a straightforward calculation we get

$$d\omega = -(g_{ij} - l_i l_j)dx^i \wedge \frac{\delta y^j}{F}. \tag{1}$$

### 1.7. Gradient vector field

For a Riemannian manifold $(SM, \widetilde{g})$, the gradient vector field of a function $f \in C^\infty(SM)$ is given by $\widetilde{g}(\nabla f, \widetilde{X}) = df(\widetilde{X}), \forall \widetilde{X} \in \chi(SM)$. Using the local coordinate system $(x^i,[y^i])$ for $SM$, the vector field $\widetilde{X} \in \chi(SM)$ is given by $\widetilde{X} = X^i(x,y)\frac{\delta}{\delta x^i} + Y^i(x,y)F\frac{\partial}{\partial y^i}$ where $X^i(x,y)$ and $Y^i(x,y)$ are $C^\infty$ functions on $SM$. A simple calculation shows that locally

$$\nabla f = g^{ij} \frac{\delta f}{\delta x^i} \frac{\delta}{\delta x^j} + F^2 g^{ij} \frac{\partial f}{\partial y^i} \frac{\partial}{\partial y^j}.$$

The norm of $\nabla f$ with respect to the Riemannian metric $\widetilde{g}$ is given by

$$|\nabla f|^2 = \widetilde{g}(\nabla f, \nabla f) = g^{ij} \frac{\delta f}{\delta x^i} \frac{\delta f}{\delta x^j} + F^2 g^{ij} \frac{\partial f}{\partial y^i} \frac{\partial f}{\partial y^j}. \tag{2}$$

## 2. EXTENSION OF SOME DEFINITIONS TO FINSLER MANIFOLDS

In what follows, $(M,g)$ denotes a connected Finsler manifold of class $C^1$ with dimension $n \geq 2$. Let $(SM, \widetilde{g})$ be its Riemannian Sphere bundle.
We consider the volume element $\eta(g)$ on $SM$ defined as follows:

$$\eta(g) := \frac{(-1)^N}{(n-1)!} \omega \wedge (d\omega)^{n-1}, \tag{3}$$

where $N = \frac{n(n-1)}{2}$ and $\omega$ is the Hilbert form of $F$ (This volume element was used for the first time in Finsler geometry by Akbar-Zadeh in his thesis [11] and [17]). Let $C(M)$ be the linear space of continuous real valued functions on $M$, $u \in C(M)$ and $u^V$ its vertical lift on $SM$. For $M$, compact or not, we denote by $H(M)$ the set of all functions in $C(M)$, admitting a generalized $L^n$-integrable gradient $\nabla u^V$ satisfying

$$I(u,M) = \int_{SM} |\nabla u^V|^n \, \eta(g) < \infty.$$

If $M$ is non-compact let us denote by $H_0(M)$ the subspace of functions $u \in H(M)$ for which the vertical lift $u^V$ has a compact support in $SM$. A *relatively compact* subset is a subset whose closure is compact. A function $u \in C(M)$ will be called *monotone* if for any relatively compact domain $D$ of $M$

$$\sup_{x \in \partial D} u(x) = \sup_{x \in D} u(x); \quad \inf_{x \in \partial D} u(x) = \inf_{x \in D} u(x).$$

We denote by $H^*(M)$ the set of monotone functions $u \in H(M)$.

We define notion of capacity as follows:

**Definition 2.1.** Capacity of a compact subset $C$ of a non-compact Finsler manifold $M$ is defined by

$$Cap_M(C) := \inf_u I(u, M),$$

where the infimum is taken over the functions $u \in H_0(M)$ with $u = 1$ on $C$ and $0 \leq u(x) \leq 1$ for all $x$, these functions are said to be admissible for $C$.

The non-compactness condition of M is a necessary condition. In fact, if M is compact, then by putting $u = 1$ in $H_0 M$ we have $I(u, M) = 0$, therefore the capacity of all subsets is zero and there is nothing to say.

A *relative continuum* is a closed subset $C$ of $M$ such that $C \cup \{\infty\}$ is connected in Alexandrov's compactification $\overline{M} = M \cup \{\infty\}$. To avoid ambiguities, the connected closed sets of $M$ that are not reduced to one point will be called *continua*. In what follows we want to associate conformal invariant function, which is determined entirely by the conformal structure of manifold $M$, at every double point of $M$.

**Definition 2.2.** Let $(M, g)$ be a Finsler manifold. For all $(x_1, x_2)$ in $M^2 := M \times M$ we set

$$\mu_M(x_1, x_2) = \inf_{C \in \alpha(x_1, x_2)} Cap_M(C),$$

where $\alpha(x_1, x_2)$ is the set of all compact continua subsets of $M$ containing $x_1$ and $x_2$.

## 3. CONFORMAL PROPERTY OF CAPACITY

**Lemma 3.1.** Let $(M, g)$ and $(M', g')$ be two conformal related Finsler manifolds, then there exist an orientation preserving diffeomorphism between their sphere bundles.

**Proof:** Let $f : (M, g) \longrightarrow (M', g')$ be a diffeomorphism between two Finsler manifolds. We define a mapping $h$ between their sphere bundles as follows $h : SM \longrightarrow SM'$, where $h(x, [y]) = (f(x), [f_*(y)])$, and $f_*$ is the differential map of $f$. Since $f_*$ is a linear map, $h$ is well defined. If $f$ is conformal then $f^* g' = \lambda g$, where $\lambda$ is a positive real valued function on $M$ and for components of Finsler metrics $g$ and $g'$ defined on $TM$ and $TM'$ we have $\lambda g = f^* g' = f^*(g'_{ij} dx'^i dx'^j)$, by definition $(f_*)^* g'_{ij}(f^* dx'^i)(f^* dx'^j) = (f_*)^* g'_{ij} dx^i dx^j$, and therefore $(f_*)^* g'_{ij} = \lambda g_{ij}$ or equivalently, $h^* g'_{ij} = \lambda g_{ij}$. Let $\omega'$ be the Hilbert form related to the Finsler metric $g'$. By definition

$$\omega' = g'_{ij} \frac{y'^j}{F} dx'^i = g'_{ij} \frac{y'^j}{\sqrt{g'_{mn} y'^m y'^n}} dx'^i.$$

Therefore,

$$h^* \omega' = h^*(g'_{ij}) \frac{h^*(y'^j)}{\sqrt{h^*(g'_{mn} y'^m y'^n)}} h^*(dx'^i) = \sqrt{\lambda} \omega. \tag{4}$$

By applying $h^*$ to (1) we get by straight forward calculation

$$h^* d\omega' = \sqrt{\lambda} d\omega. \tag{5}$$

So if $\eta(g)$ and $\eta(g')$ denote the volume elements of $SM$ and $SM'$ respectively, then from (3), (4) and

(5) we get

$$h^*(\eta(g')) = (\sqrt{\lambda})^n \eta(g). \tag{6}$$

Therefore $h$ is an orientation preserving diffeomorphism.

**Lemma 3.2.** Let $f$ be a diffeomorphism between Finsler manifolds $(M, g)$ and $(M', g')$, and $h$ a mapping between their sphere bundles with Sasaki metrics, $(SM, \tilde{g})$ and $(SM', \tilde{g}')$. If $u \in H_0(M')$ then we have

1. $|\nabla u^V|^n = (g'^{ij} \frac{\delta u^V}{\delta x'^i} \frac{\delta u^V}{\delta x'^j})^{\frac{n}{2}}$,
2. $(u \circ f)^V = u^V \circ h,$
3. $h^* \frac{\delta u^V}{\delta x'^i} = \frac{\delta (u \circ f)^V}{\delta x^i}.$

Therefore, the following diagram is commutative:

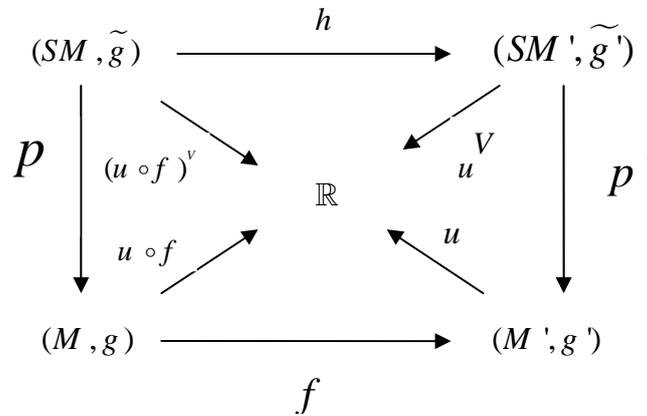

Diagram 1.

**Proof:**
1. Since the vertical lift of $u \in H_0(M')$ is a function of position alone, $\frac{\partial u^V}{\partial y^i} = 0$. Therefore the first assertion follows from (2).
2. Let's consider the projections $p : SM \to M$ and $p' : SM' \to M'$. The vertical lifts of $u$ and $u \circ f$, are by definition, $u^V(x', [y']) = u \circ p'(x', [y']) = u(x')$ and

$$(u \circ f)^V (x, [y]) = (u \circ f) \circ p(x, [y]) = (u \circ f)(x).$$

From which we have

$$(u \circ f)^V (x, [y]) = (u \circ f)(x) =$$
$$u^V (f(x), [f_*(x)]) = u^V (h(x, [y])) = u^V \circ h(x, [y]).$$

This proves the assertion $(2)$.

3. By definition of $h^*$ we have $h^*(\frac{\delta}{\delta x'^i} u^V) = h^*(\frac{\delta}{\delta x'^i}).h^*u^V = \frac{\delta}{\delta x^i}.(u^V \circ h),$ and from (2) we get assertion $(3)$.

Now we are in a position to prove the following theorem:

**Theorem 3.3.** Let $(M, g)$ be a connected non-compact Finsler manifold, then the capacity of a compact set on $M$ is a conformal invariant.

**Proof:** We show that the notion of capacity depends only on the conformal structure of $M$, or equivalently, for any conformal map $f$ from Finsler manifold $(M, g)$ onto another Finsler manifold $(M', g')$, we have

$$Cap_M(C) = Cap_{M'}(f(C)).$$

Since $SM$ and $SM'$ are two smooth, orientable manifolds with boundary, then for a smooth, orientation preserving diffeomorphism function $h: SM \longrightarrow SM'$ defined in Lemma 3.1, clearly (by a classical result in differential Geometry, [18]) we have

$$\int_{SM'} \omega = \int_{SM} h^*\omega, \qquad \omega \in \Omega^{2n-1} SM'.$$

So we get,

$$I(u, M') = \int_{S(M')} |\nabla u^V|^n \eta(g') = \int_{SM} h^*(|\nabla u^V|^n \eta(g')). \qquad (7)$$

Using Lemma 3.2, a straightforward calculation shows that

$$h^* |\nabla u^V|^n = (\sqrt{\lambda})^{-n} |\nabla (u \circ f)^V|^n. \qquad (8)$$

Using (6) in Lemma 3.1, and relations (7) and (8) we get

$$I(u, M') = \int_{SM} |\nabla (u \circ f)^V|^n \eta(g) = I(u \circ f, M). \qquad (9)$$

Let $C$ be a compact set in $M$, then we have

$$Cap_M(C) = \inf_{v \in H_0 M, v|_C = 1} I(v, M), Cap_{M'}(f(C)) = \inf_{u \in H_0 M', u|_{f(C)} = 1} I(u, M').$$

Put

$$A = \{I(v, M) \mid v \in H_0 M, v|_C = 1\},$$

$$B = \{I(u, M') \mid u \in H_0 M', u|_{f(C)} = 1\}.$$

We first show that $B \subseteq A$. For all $I(u, M') \in B$, we easily have the following assertions.
- Since $support(u^V)$ is compact in $SM'$, $h^{-1}(support(u^V)) = support(u \circ f)^V$ is compact in $SM$ and by definition $u \circ f \in H_0(M)$.
- $(u \circ f)|_C = 1$ since $u|_{f(C)} = 1$.
- From (9) we have $I(u \circ f, M) = I(u, M')$.

Therefore, $I(u \circ f, M) \in A$ and $B \subseteq A$. By the same argument we have $A \subseteq B$. Hence, $Cap_M(C) = Cap_{M'}(f(C))$.

Theorem 3.3, implies that the function $\mu_M$ is invariant under any conformal mapping. More precisely, if $f$ is a conformal mapping between Finsler manifolds $(M, g)$ and $(M', g')$, then for all $x_1, x_2 \in M$ we have

$$\mu_M(x_1, x_2) = \mu_{M'}(f(x_1), f(x_2)),$$

In the Riemannian geometry this function is of general interest in the study of global conformal geometry, which can be the subject of further studies in Finsler geometry.